\documentclass[12pt] {article}
\usepackage[cp1251]{inputenc}
\usepackage[english, russian] {babel}
\usepackage {amsmath}
\usepackage {amsthm, amssymb}

\usepackage[dvips]{graphicx}
\usepackage {graphicx}
\newcommand{\aff}[1]{\mathop{aff}\left(#1\right)}
\newcommand{\lin}[1]{\mathop{lin}\left(#1\right)}
\newcommand{\notparallel}{\not\,\parallel}
\newcommand{\epi}{\mathop{epi}}

\renewcommand{\int}[1]{\mathop{Int}\left(#1\right)}

\newcommand{\p}{\mathbf{p}}
\newcommand{\m}{\mathbf{m}}

\newcommand{\mcP}{\mathcal{P}}
\newcommand{\mcM}{\mathcal{M}}
\newcommand{\mcT}{\mathcal{T}}
\newcommand{\mcG}{\mathcal{G}}
\newcommand{\mcA}{\mathcal{A}}

\newcommand{\parenth}[1]{\left(#1\right)}
\newcommand{\brackets}[1]{\left[#1\right]}
\newcommand{\braces}[1]{\left\{#1\right\}}

\newcommand{\column}[2]{\left({#1 \atop #2}\right)}

\newtheorem{theorem}{Теорема}[section]
\newtheorem*{classicTheorem}{Теорема}

\theoremstyle{plain}

\newtheorem{proposition}[theorem]{Предложение}
\newtheorem{lemma}[theorem]{Лемма}
\newtheorem{corollary}[theorem]{Следствие}
\newtheorem{definition}[theorem]{Определение}

\newtheorem*{remark}{Замечание}

\begin{document}

\title{Геометрия подъёмов разбиений евклидовых пространств}
\author{А.А.~Гаврилюк}


\date{\today}

\begin{abstract}
  В данной работе даётся полное обоснование метода канонических нормировок полиэдральных комплексов
  для подъёма разбиений евклидова пространства. Предложен новый комбинаторно-геометрический подход
  к построению женератрисы разбиения на основе операции подъёма грани до уже построенного соседа.
  
  На основе изложенного общего подхода даётся новое, существенно сокращённое по сравнению с 
  имеющимися, геометрическое доказательство фундаментальной теоремы из теории параллелоэдров, что 
  для данного параллелоэдра $P$ гипотеза Вороного верна тогда и только тогда, когда для 
  соответствующего разбиения $\mathcal T_P$ существует каноническая нормировка.
  
  \bigskip
  \noindent {\bf Ключевые слова:} параллелоэдры, гипотеза Вороного, женератриса, каноническая 
    нормировка
  
  \bigskip
  \noindent УДК 514.174 + 514.87
\end{abstract}

\maketitle

\section{Введение, основные понятия}
\subsection{Разбиение на параллелоэдры и подъёмы разбиений}

  Основы теории параллелоэдров были заложены в работе кристаллографа Е.С. Фёдорова \cite{Fedorov} 
  (1885) и были сформулированы для размерности 3. В общем случае $d$-мерным {\it параллелоэдром} 
  называется выпуклый евклидов многогранник $P$, который своими параллельными копиями разбивает 
  пространство ${\mathbb R^d}$ нормальным образом, то есть так, что пересечение любых двух 
  многогранников либо пусто, либо является их общей целой гранью некоторой размерности. Важное 
  свойство, связывающее параллелоэдры с теорией чисел --- то, что центры параллелоэдров образуют 
  $d$-мерную целочисленную решётку $\Lambda^d(P)$ в некотором базисе пространства $\mathbb E^d$ 
  \cite{Dolbilin_new}.

  Г. Минковский \cite{Minkovski} доказал замечательную теорему, дающую необходимые условия того, 
  что выпуклый многогранник является параллелоэдром:

  \begin{classicTheorem}[Минковский]\label{theoremMinkovski}~ {\it 
	  Если $P$ --- $d$-мерный параллелоэдр, то \quad (1) многогранник $P$ центрально-симметричен, 
	  \quad (2) все его гиперграни центрально-симметричны.
  } \end{classicTheorem}

  \noindent Позднее эти условия были дополнены ещё одним:

  \begin{classicTheorem}[Венков, Делоне]\label{condDelone}~ {
    Проекция $d$-мерного параллелоэдра $P$ вдоль каждой его $(d-2)$-грани на дополнительную 
    2-мерную плоскость является либо параллелограммом, либо центрально-симметричным 
    шестиугольником.
  } \end{classicTheorem}
 
  Долгое время это условие также приписывалось Минковскому, однако ни в оригинальной работе 
  \cite{Minkovski}, ни в прочих работах Минковского этот результат не только не сформулирован явно, 
  но, судя по всему, и вовсе не упоминается. Впервые это условие упоминается, по всей видимости, в 
  работе Б.Н. Делоне \cite{Delone}, однако как отдельная теорема оно также не формулируется. Лишь в
  статье Б.А. Венкова \cite{Venkov} в 1954 году (то есть через 25 лет после публикации работы 
  Делоне и более чем через 50 лет после выхода в свет работы Минковского) эти условия были собраны 
  вместе и было показано, что выполнения этих трёх условий достаточно, чтобы выпуклый многогранник 
  являлся параллелоэдром.

  Одна из основных задач теории параллелоэдров --- нахождение алгоритма, перечисляющего для данной 
  размерности все комбинаторные типы параллелоэдров. Эта задача до сих пор остаётся нерешённой. 
  Г.Х. Вороной \cite{Voronoy} построил теорию параллелоэдров специального вида, ныне носящих его 
  имя, и привёл алгоритм перечисления всех их комбинаторных типов. Под \emph{параллелоэдром 
  Вороного} $P_V$ мы понимаем область Дирихле некоторой точки $\mathbf{p}$ из $d$-мерной целочисленной 
  решётки $\Lambda^d \subset \mathbb R^d$:
  $$P_V = \left\{\mathbf{x} \in \mathbb E^d : \|\mathbf{x} - \mathbf{p}\| \leqslant \|\mathbf{x} - 
    \mathbf{q}\| \; \forall \mathbf{q} \in \Lambda^d \right\}$$
  Легко проверяется, что такая область действительно является параллелоэдром. Вороной высказал 
  гипотезу о том, что любой параллелоэдр общего вида аффинно-эквивалентен некоторому параллелоэдру 
  $P_V$. Доказательство этой гипотезы будет означать, что алгоритм Вороного перечисляет и все
  комбинаторные типы параллелоэдров в общем случае.

  Нормальное разбиение $\mathcal T_P$ называется \emph{примитивным} в $n$-мерной грани $F^n$, 
  $0 \leqslant n \leqslant d-1$, если в этой грани сходятся в точности $(d+1-n)$ $d$-мерных ячеек 
  разбиения, что является минимально возможным количеством для заданных $n$ и $d$. Саму грань $F^n$
  в таком случае мы будем называть \emph{\it примитивной}. Если многогранник примитивен во всех 
  своих $k$-гранях, то он называется \emph{$k$-примитивным}.

  Вороной доказал свою гипотезу для 0-примитивных параллелоэдров (называются просто 
  \emph{примитивными}). Для разбиения на такие параллелоэдры он построил специальную функцию --- 
  женератрису, график $\mathcal G_P$ которой является $d$-мерной границей неограниченного 
  $(d+1)$-мерного полиэдра и проектируется в точности в разбиение $\mathcal T_P$. $\mathcal G_P$ 
  при этом естественно рассматривать как подъём разбиения $\mathcal T_P$ над $\mathbb E^d$. 
  
  О.К. Житомирский \cite{Zhitomirsky} усилил теорему Вороного и доказал гипотезу для очень широкого
  класса --- $(d-2)$-примитивных параллелоэдров. Это ограничение эквивалентно тому, что в каждой 
  $(d-2)$-мерной грани сходятся ровно 3 параллелоэдра и является, в некотором смысле, самым слабым 
  ограничением (то есть дающим самый сильный результат), определяемым через примитивность. Это 
  следует из простого наблюдения, что всякий $k$-примитивный параллелоэдр $(0\leqslant k\leqslant 
    d-2)$, также является $(k+1)$-примитивным.

  Много позднее, более детальным рассмотрением локальной структуры были получены различные усиления 
  случая Житомирского: для \emph{3-неразложимых} параллелоэдров \cite{Ordin}, при {\it односвязной 
  $\delta$-поверхности} \cite{GarGavMag}, для суммы параллелоэдра и отрезка \cite{Grishukhin, 
  Magazinov}. Альтернативный подход к доказательству гипотезы предложил Р. Эрдал \cite{Erdahl} для 
  параллелоэдров, являющихся зоноэдрами.

  Конструкция подъёма разбиений рассматривается и для общих локально-конечных разбиений 
  евклидова пространства. В случае примитивных разбиений подъём можно совершать без дополнительных 
  конструкций. Для этого сперва надо сделать произвольный согласованный подъём двух начальных 
  смежных $d$-ячеек. Далее достаточно добавлять к уже поднятому участку новые подъёмы, смежные 
  с какой-то одной поднятой ячейкой по гиперграни, а с другой --- хотя бы по ребру. Как показано в
  \cite{Davis}, при $d \geqslant 3$ для конечных примитивных разбиений такое построение всегда 
  корректно определяет однозначный подъём.
  
  Для общих разбиений такое построение может оказаться несогласованным, поэтому требуются 
  дополнительные ограничения. Таким ограничением обычно служит существование 
  \emph{ортогонально-дуального} множества: каждой ячейке разбиения сопоставляется точка такого 
  множества. При этом отрезок, соединяющий точки для смежных ячеек, должен быть ортогонален их 
  общей гиперграни \cite{McMullen}. На практике построение такого множества (и в такой 
  формулировке) бывает затруднительно.
  
  Канонические нормировки, впервые появившиеся в ``чистом'' виде, по всей видимости, в 
  \cite{Ordin}, с одной стороны эквивалентны конструкции ортогонально-дуальных множеств, с другой 
  --- позволяют более эффективно работать с локальной структурой разбиений. Не смотря на 
  появившиеся в последнее время результаты, полученные при помощи этой техники \cite{Ordin, 
  GarGavMag}, полного её обоснования никогда не было приведено: имеющиеся изложения используют
  другие эквивалентные конструкции и исключительно аналитическое их описание.
  
  В данной работе приведено полное самодостаточное обоснование метода канонических нормировок (с
  единственной существенно используемой ссылкой на теорему \ref{theoremQualityTr}), приведён новый
  геометрический подход к построению женератрисы разбиения общего вида на основе операции подъёма
  грани до ранее поднятого соседа (определение \ref{definitionLift}). Приведены геометрические 
  доказательства (теорема \ref{theoremManifoldGeneratris}, теорема \ref{theoremGeneratrissConvex})
  того, что в результате такого построения получается $d$-мерное полиэдральное многообразие, 
  ограничивающее выпуклый полиэдр. Далее на основе этого построения приводится существенно более
  короткое, чем другие аналогичные (\cite{DezaGrishukhin}, в частном случае \cite{Voronoy}) 
  геометрически прозрачное доказательство фундаментальной теоремы \ref{theoremMainEquivalence}
  о том, что для параллелоэдра $P$ верна гипотеза Вороного тогда и только тогда, когда для 
  соответствующего разбиения $\mathcal T_P$ существует каноническая нормировка. Применение подхода
  иллюстрируется на доказательстве теоремы Житомирского.

\subsection{Функции приращения}

	Рассмотрим некоторый однородный $n$-мерный полиэдральный комплекс $\mathcal{K}^n$. \emph{Цепью} в 
	$\mathcal{K}^n$ называется произвольная последовательность $n$-мерных ячеек $[c_1,\ldots,c_k]$, в 
	которой любые две последовательные ячейки \emph{\it смежны}, то есть имеют общую $(n-1)$-грань. 
	\emph{(Комбинаторным) циклом} называется произвольная цепь, у которой первая и последняя ячейки 
	смежны. Инцидентность ячеек комплекса комплекса тут и далее мы индуктивно определяем через 
	смежность по $(d-1)$-граням:
	
	\begin{definition}[Локальная структура комплекса]\label{definitionComplex}~ 
		\begin{itemize}
			\item Если $n$-мерные ячейки $c$ и $c'$ смежны, то есть имеют общую $(d-1)$-грань $F^{d-1}$, то 
				каждая подгрань $G \subseteq F^{d-1}$ считается общей для $c$ и $c'$
			\item Если для грани $G$ и ячеек $c$ и $c'$ найдётся соединяющая их цепь 
				$[c = c_1, c_2,\ldots, c_k = c']$ такая, что $G \subseteq c_1 \bigcap c_2$, а значит является
				гранью и $c$, и $c'$, $G \subseteq c_2 \bigcap c_3$, а значит является гранью и $c_3$, и так 
				далее до $G \subseteq c_{k-1} \bigcap c_k$, а значит является гранью и $c' = c_k$, то $G$
				считается общей гранью для $c$ и $c'$.
		\end{itemize}
	\end{definition}
	Такая цепь $n$-мерных ячеек $[c_1, c_2,\ldots, c_k]$, все из которых содержат некоторую общую 
	$k$-мерную грань $G$, называется \emph{$k$-примитивной}.
	
	\begin{remark} В случае разбиений пространств такое определение инцидентности совпадает со 
	  стандартным теоретико-множественным определением через пересечение многогранников как точечных 
	  множеств. Однако для более изощрённых структур, например подмножеств
		разбиений пространства или абстрактно заданных полиэдральных многообразий не вложенных в 
		евклидово пространство, теоретико-множественное определение может оказаться неудобным или
		вовсе неопределённым. Описанное же определение работает всегда, когда заданы смежности по
		$(d-1)$-граням (см., например, \cite{Alexandrov}).
	\end{remark}

  \emph{Нормировкой} комплекса $\mathcal{K}^n$ будем называть произвольную функцию $s$, определённую 
  на $\mathcal F^n = \mathcal F^n_{\mathcal K^n}$ --- множестве $n$-мерных граней $\mathcal K^n$.
	
	\begin{definition}[Мультипликативные приращения]\label{definitionScaling}  
		Пусть для каждой упорядоченной пары смежных $n$-мерных ячеек $c_1$ и $c_2$ определена функция 
		$g:~[c_1, c_2]~\to~\mathbb R_+$ такая, что $g[c_1, c_2] \cdot g[c_2, c_1] = 1$. 
		\emph{Мультипликативным приращением} вдоль цепи $[c_1, c_2,\ldots, c_k]$ называется продолжение 
		функции $g$ по естественному правилу 
		$$g[c_1, c_2,\ldots, c_k] = g[c_1, c_2]g[c_2, c_3]\cdot\ldots\cdot g[c_{k-1}, c_k]$$
	\end{definition}
	Если зафиксировать некоторые ячейку $c_0$ и положительное значение $s_0$, то правило
	$s(c) := s_0\cdot g[c_0, c_1,\ldots, c]$ по мультипликативному приращению задаёт на 
	$\mathcal{K}^n$ положительно определённую нормировку, которая, вообще говоря, зависит от цепи, 
	соединяющей ячейки $c_0$ и $c$.
	\begin{definition}[Аддитивные приращения]\label{definitionAddScaling}  
		\emph{Аддитивным приращением} называется функция $g:~[c_1, c_2]~\to~\mathbb R^k, 
		  k\in \mathbb N$ с условием $g[c_1, c_2] + g[c_2, c_1] = 0$, продолженная на цепи комплекса 
		  по правилу $g[c_1, c_2,\ldots, c_k] := g[c_1, c_2]+\ldots+g[c_{k-1}, c_k]$.
	\end{definition}
	Аналогично, для фиксированной ячейки $c_0$ и произвольного значения $v_0\in \mathbb R^k$ правило
	$s(c) := v_0 + g[c_0, c_1,\ldots, c]$ задаёт вектор-значную нормировку.
	
	Очевидно, что нормировка определена приращениями однозначно с точностью до 
	начального значения тогда и только тогда, когда приращения вдоль любых двух цепей с общими 
	началом и концом, равны. Проверку этого свойства существенно упрощает теорема о переносе свойства 
	С. Рышкова и К. Рыбникова \cite{RyshRyb}. Авторы доказали свою теорему для более широкого класса 
	комплексов, чем потребует данная работа, поэтому мы приведём упрощённую версию:  

	\begin{theorem}[О переносе свойства]\label{theoremQualityTr}~{ 
		Рассмотрим $n$-мерный комплекс $\mathcal K^n$ в $\mathbb{R}^d$, который является либо 
		локально-конечным разбиением $\mathbb{R}^d$, либо остовом размерности $n \geqslant 2$ такого 
		разбиения. Пусть на всех парах смежных $n$-мерных ячеек в $\mathcal K^n$ заданы приращения
		$g[c_1, c_2]$ и по некоторой начальной ячейке $c_0$ и значению на ней $s_0$ задана нормировка
		$s$ данного комплекса. Функция $s$, как функция только от ячеек разбиения, существует (и 
		определена однозначно с точностью до выбора начального значения) тогда и только тогда, 
		когда приращение вдоль всех $(n-2)$-примитивных циклов равно $1$ в случае мультипликативных
		приращений, либо $0$ в случае аддитивных приращений.}
	\end{theorem}

\subsection{Канонические нормировки}

	Рассмотрим нормальное (то есть грань-в-грань) разбиение $\mathcal T$ пространства $\mathbb R^d$ на 
	строго выпуклые многогранники. Пусть задана положительная нормировка $s:~{\mathcal F}^{d-1}~\to~
	  \mathbb R_+$ $(d-1)$-мерного остова $\mathcal F^{d-1}$ разбиения $\mathcal T$. Пусть 
	$F_1, F_2, \ldots, F_k$ --- все гиперграни разбиения, которые содержат $(d-2)$-грань $F^{d-2} 
	  \subset \mathcal F^{d-2}$. На них естественным образом определены два противоположных 
	направления обхода, которые определяются проекцией гиперграней $F_i$ и самой $F^{d-2}$ в двумерную 
	плоскость дополнительную к аффинной оболочке $\aff{F^{d-2}}$. Будем считать, что нумерация 
	$F_1, F_2, \ldots, F_k$ задана одним из этих двух направлений обхода, и, согласованно с тем же
	направлением, выбраны единичные нормали $\mathbf{n}_1, \ldots, \mathbf{n}_k$ к	данным гиперграням. 
	{\it Кручением} $\Delta_s$ нормировки $s$ вокруг $(d-2)$-грани $F^{d-2}$ назовём величину 
	$$\Delta_s(F^{d-2}) := \sum^k_{i=0} s(F_i)\mathbf{n}_i$$
	Из определения следует, что кручение определено с точностью до знака и, вообще говоря, зависит от
	выбора направления обхода. Однако нас интересует лишь равенство или неравенство кручения нулю, что
	от направления обхода не зависит. 

	Нормировка $s$ множества гиперграней $\mathcal F^{d-1}$ разбиения $\mathcal T$ называется 
	\emph{канонической}, если для любой $(d-2)$-грани $F^{d-2} \subset \mathcal T$ выполнено 
	$\Delta_s(F^{d-2}) = 0$. 

	{\it Локальной} канонической нормировкой $s: {\mathcal S}^{d-1} \to \mathbb R_+$ подкомплекса 
	${\mathcal S^{d-1}} \subset \mathcal T$ мы называем нормировку с нулевым кручением для всякой 
	$(d-2)$-грани $F^{d-2} \subset {\mathcal S}^{d-1}$, звезда $St_{\mathcal T}(F^{d-2})$ которой, 
	определённая по комплексу $\mathcal T$, принадлежит и подкомплексу	$\mathcal S^{d-1}$.

	Следующая лемма --- одно из ключевых утверждений при рассмотрении канонических нормировок:
	\begin{lemma}\label{lemmaUniqTri}~{
		Пусть $F^{d-2}$ --- $(d-2)$-грань разбиения $\mathcal T$, в которой сходятся ровно 3 гиперграни 
		$F_1, F_2$ и $F_3$ из $\mathcal T$. Пусть $\mathcal S$ и $\mathcal S'$ --- два подкомплекса в 
		$\mathcal T$, оба целиком содержат звезду $St_{\mathcal T}(F^{d-2})$. Пусть на этих	подкомплексах 
		заданы (локальные) канонические нормировки $s$ и $s'$. Тогда выполнено равенство:
		$$\frac{s'(F_1)}{s'(F_2)} = \frac{s(F_1)}{s(F_2)}$$
	}
	\end{lemma}
	
	\smallskip
	\noindent {\it Доказательство.}
	По определению канонической нормировки имеем: 
	$$s(F_1)\mathbf{n}_1 + s(F_2)\mathbf{n}_2 + s(F_3)\mathbf{n}_3 = 
		s'(F_1)\mathbf{n}_1 + s'(F_2)\mathbf{n}_2 + s'(F_3)\mathbf{n}_3 = 0, ~s(F_i) > 0, s'(F_i) > 0$$
	Отметим, что из строгой выпуклости ячеек $\mathcal T$ следует, что среди $F_1, F_2, F_3$ не может 
	быть двух параллельных гиперграней. Таким образом, $\mathbf{n}_1$ и $\mathbf{n}_3$ не параллельны.
	Значит, если $\mathbf{n}_2$ раскладывается по этим векторам, то единственным образом. В частности, 
	коэффициент при $\mathbf{n}_1$ в таком разложении определён однозначно. С другой стороны, из 
	выписанных равенств имеем: 
	$$\mathbf{n}_2 = - \frac{s(F_1)}{s(F_2)}\mathbf{n}_1 - \frac{s(F_3)}{s(F_2)}\mathbf{n}_3
		= - \frac{s'(F_1)}{s'(F_2)}\mathbf{n}_1 - \frac{s'(F_3)}{s'(F_2)}\mathbf{n}_3$$
	\hfill $\Box$ \par\bigskip

	Вернёмся от общих разбиений пространства к разбиению на копии заданного параллелоэдра $P$. Для
	таких разбиений кроме примитивных граней следует выделить важный класс стандартных граней, 
	введённых в \cite{Dolbilin_new}. Грань $F^n$ разбиения $\mathcal T_P$ называется 
	\emph{стандартной}, если её можно представить в виде $F^n = P \bigcap P'$ для некоторых ячеек $P$ 
	и $P'$ разбиения. Из центральной симметрии параллелоэдров следует, в частности, что стандартные 
	грани центрально симметричны, и их центры симметрии являются также центрами симметрии разбиения 
	$\mathcal T_P$ и решётки $\Lambda^d(P)$.

	Из условия Венкова-Делоне следует следующая простая классификация схождений параллелоэдров в 
	$(d-2)$-гранях: 

	\begin{proposition}\label{propositionDminus2}~{
		Произвольная $(d-2)$-грань разбиения $\mathcal T_P$ является либо примитивной, в которой сходятся 
		3 гиперграни разбиения, либо стандартной, в которой сходятся 4 гиперграни $F_1, F_2, F_3, F_4$,
		причём $F_1 || F_3$ и $F_2 || F_4$.
	}
	\end{proposition}

	Гиперграни параллелоэдров являются частным случаем стандартных граней. То есть центр каждой 
	гиперграни является также центром симметрии разбиения. Рассмотрим две смежные по $(d-2)$-грани
	$F^{d-2}_{12}$ гиперграни $F_1, F_2$ некоторого параллелоэдра $P_1 \subset 
	  \mathcal T_P$. Пусть $F_1' \subset \mathcal T_P$ --- это образ $F_1$ при 
	симметрии относительно центра $F_2$. Отметим, что $F_1'$ и $F_2$ также смежны и принадлежат 
	параллелоэдру $P'_1 = \mathop{Sym}_{c(F_2)}\parenth{P_1}$ данного разбиения. Будем называть такую
	пару гиперграней $F_1$ и $F'_1$ \emph{накрест лежащими} относительно (смежной с ними обеими) 
	$F_2$.
	
	\begin{lemma}\label{lemmaCross}~{ 
	  Пусть задана каноническая нормировка $s$ разбиения $\mathcal T_P$, и пусть $\brackets{F_1, F_2, 
	    \ldots, F_6}$ --- произвольный 6-поясок параллелоэдра $P_0 \subset \mathcal T_P$
	  Пусть $F_0$ --- третья гипергрань в $\mathcal T_P$, содержащая $(d-2)$-грань 
	  $F_2 \bigcap F_3$ (помимо самих $F_2$ и $F_3$). Тогда \quad (1) нормировки на противоположных 
	  гипергранях из 6-пояска равны, \quad (2) нормировки на накрест лежащих гипергранях $F_1$ и $F_0$ 
	  равны.
	} \end{lemma}
	
	\smallskip
	\noindent {\it Доказательство.}
	  Докажем сначала пункт 2. $F_1$ и $F_0$ --- накрест лежащие относительно $F_2$. Применим к 
	  разбиению $\mathcal T_P$ симметрию $c\parenth{F_2}$. $F_1$ и $F_0$ поменяются местами, сама 
	  $F_2$ и $\mathcal T_P$ перейдут в себя. Нормировка $s$ перейдёт в некоторую нормировку $s'$
	  разбиения $\mathcal T_P$, причём $s'(F) = s\parenth{\mathop{Sym}_{c\parenth{F_2}}F}$. 
	  Из леммы \ref{lemmaUniqTri} и симметрии имеем $\frac{s(F_1)}{s(F_2)} = \frac{s'(F_1)}{s'(F_2)} =
	    \frac{s(F_0)}{s(F_2)}$. Значит $s(F_1) = s(F_0)$ и пункт 2 доказан. Согласно этому пункту 
	  также верно $s(F_4) = s(F_0)$. Значит $s(F_4) = s(F_1)$ и пункт 1 также доказан.
	\hfill $\Box$ \par\bigskip
	
	Из леммы \ref{lemmaUniqTri} и предложения \ref{propositionDminus2} следует, что канонические 
	нормировки звезды $St_{\mathcal T_P}(F^{d-2})$ произвольной $(d-2)$-грани $F^{d-2}$ разбиения 
	$\mathcal T_P$ описываются следующим образом:

	\begin{proposition}\label{propositionDminus2Scaling}~{
		Если $F^{d-2}$ примитивна, то локальная каноническая нормировка $s$ её звезды единственна с 
		точностью до положительного общего множителя, иначе $F^{d-2}$ стандартна и нормировка $s$ четырёх
		содержащих её гиперграней $F_1 || F_3$ и $F_2 || F_4$ является канонической тогда и только тогда,
		когда выполняются равенства 
		$$s(F_1) = s(F_3) > 0, ~ s(F_2) = s(F_4) > 0$$
	}
	\end{proposition}

	Далее нам потребуется усиление этого утверждения:
	\begin{theorem}[Делоне \cite{Delone}]\label{theoremD_3Scaling}~{
	  Звёзды произвольных $(d-2)$- и $(d-3)$-граней параллелоэдров имеют канонические нормировки (не 
	  менее одной).
	} \end{theorem}

\subsection{Склейка локальных нормировок}

	Одна из основных функций локальных канонических нормировок --- построение по ним канонических
	нормировок всего разбиения. Такое построение происходит при помощи \emph{склейки} нормировок. 
	Пусть каноническая нормировка $s$ задана на комплексе $\mathcal K$, каноническая нормировка $s'$ 
	--- на комплексе $\mathcal K'$. Очевидно, что на $\mathcal K'$ также определены канонические 
	нормировки вида $ks'$ для произвольного $k>0$. Для склейки $s$ и $s'$ подбирается такое $k$,
	чтобы 
	\begin{itemize}
		\item $s$ и $ks'$ совпадали на ${\mathcal K}\bigcap{\mathcal K'}$
		\item условие каноничности для расширенной нормировки $\overline{s} = s \bigcup ks'$, то есть 
		равенство нулю кручения $\Delta_{\overline{s}}$ нормировки $\overline{s}$, выполнялось вокруг 
		произвольной $(d-2)$-грани объединённого комплекса ${\mathcal K}\bigcup{\mathcal K'}$
	\end{itemize}
	Очевидно, что процесс склейки в такой форме осуществим только если $s$ и $s'$ пропорциональны на
	${\mathcal K}\bigcap{\mathcal K'}$, что, вообще говоря, не гарантируется.
	
	Иногда удаётся упростить склейку, перейдя от нормировок к соответствующим функциям приращения:
	$g[F_1, F_2] = \frac{s\left(F_1\right)}{s(F_2)}$ для гиперграней $F_1$ и $F_2$, принадлежащих общей
	$d$-мерной ячейке разбиения и пересекающихся по некоторой $(d-2)$ грани $F^{d-2}$. Преимущество
	такого подхода в том, что если $F^{d-2}$ примитивна, то, согласно лемме \ref{lemmaUniqTri}, 
	приращение $g[F_1, F_2]$ определено однозначно и, соответственно, совпадает для нормировок $s$ и
	$s'$. В других случаях может потребоваться нетривиальное обоснование того, что $g$ и $g'$ можно 
	определить согласованно. Затем по общей функции приращения $\overline{g}$ восстанавливается общая
	нормировка $\overline{s}$, и требуется доказать её однозначность --- например, при помощи теоремы
	о переносе свойства \ref{theoremQualityTr}. Проиллюстрируем этот подход на важном случае 
	Житомирского.

\subsection{Глобальная каноническая нормировка, случай Житомирского}

	\begin{lemma}\label{lemmaCSZhitomirskii}~{
	  Для произвольного $(d-2)$-примитивного параллелоэдра $P$ существует каноническая нормировка 
	  соответствующего разбиения $\mathcal T_P$.
	} \end{lemma}

	\smallskip
	\noindent {\it Доказательство.} В случае $d = 2$ $(d-2)$-примитивными параллелоэдрами будут лишь 
	  центрально симметричные выпуклые шестиугольники. Согласно предложению 
	  \ref{propositionDminus2Scaling} существует каноническая нормировка $s$ звезды произвольной 
	  вершины $v$. Гиперграни этой звезды $St_{\mathcal T_P}(v)$ --- это одномерные рёбра $e_1, e_2, 
	    e_3$. Припишем каждому ребру разбиения, параллельному $e_1$, нормировку $s(e_1)$, ребру 
	  параллельному $e_2$ --- нормировку $s(e_2)$, ребру параллельному $e_3$ --- $s(e_3)$. Очевидно, 
	  что все кручения в таком случае равны, с точностью до знака, кручению вокруг $v$, то есть равны 
	  $0$. Значит построенная нормировка каноническая.

		Далее считаем $d\geqslant 3$. По условию разбиение $\mathcal T_P$ $(d-2)$-примитивно, значит если
		две гиперграни $F_1, F_2$ данного разбиения имеют общую $(d-2)$-грань, то, по предложению 
		\ref{propositionDminus2Scaling}, приращение $g[F_1, F_2]$ определено однозначно. Значит оно определено
		однозначно и для любой цепи на $(d-1)$-остове $\mathcal F^{d-1}$. Воспользуемся теоремой 
		\ref{theoremQualityTr} о переносе свойства.
		
		Необходимо показать, что приращение вдоль произвольного $(d-3)$-примитивного цикла $g[F_1,\ldots,
		  F_n, F_1]$ равно $1$. Все $(d-1)$-мерные ячейки этого цикла имеют некоторую общую $(d-3)$-грань
		$F^{d-3}$. По теореме \ref{theoremD_3Scaling} для $St_{\mathcal T_P}(F^{d-3})$ существует 
		каноническая нормировка $s_l$. Как было отмечено выше, приращение $g[F_i,F_{i+1}]$ при переходе 
		через примитивную $(d-2)$-грань определено однозначно. В частности, для данных $F_i$ имеем 
		$g[F_i,F_{i+1}] = \frac{s_l(F_{i+1})}{s_l(F_i)}$ (полагаем $F_{n+1} = F_1$). Поэтому приращение 
		вдоль рассматриваемого цикла равно
		$$g[F_1, \ldots, F_n, F_1] = \frac{s_l(F_2)}{s_l(F_1)} 
			\frac{s_l(F_3)}{s_l(F_2)}\cdot \ldots \cdot\frac{s_l(F_1)}{s_l(F_n)} = 1$$ 
	\hfill $\Box$ \par\bigskip

\section{Женератриса разбиения}
\subsection{Подъём многогранников}
	
	Далее считаем заданными локально-конечное нормальное разбиение $\mathcal T$ пространства 
	$\mathbb  E^d$ и каноническую нормировку $s: \mathcal F^{d-1} \to \mathbb  R_+$ разбиения 
	$\mathcal T$. По заданным разбиению и нормировке требуется построить $d$-мерную полиэдральную 
	поверхность $\mathcal G$ в расширенном пространстве $\mathbb E^{d+1}$ (ограничивающую 
	$(d+1)$-мерный полиэдр), грани которой при ортогональной проекции взаимно однозначно переходят 
	в грани тех же размерностей из $\mathcal T$.	
	
  Пусть $\mathbb E^d$ вложено в $\mathbb E^{d+1}$ как $d$-мерная гиперплоскость. \emph{Подъёмом}
  $k$-мерного многогранника $P^k \subset \mathbb E^d$, $k \leqslant d$ будем называть всякий 
  $k$-мерный многогранник $\widetilde{P^k}$ (не параллельный нормали к $\mathbb E^d$) 
  такой, что его ортогональная проекция $\Pr\parenth{\widetilde{P^k}}$ на $\mathbb E^d$ совпадает 
  с $P^k$. Тут и далее, когда это не вызывает путаницы, под $\mathbb E^d$ мы понимаем исходное 
  подпространство в $\mathbb E^{d+1}$, в котором находится поднимаемый комплекс.
  
	\begin{definition}\label{definitionLift}~{
    Пусть заданы два выпуклых многогранника $P^l_1, P^m_2 \subset \mathbb E^d$, $l, m \leqslant d$, 
    пересечение которых является гранью каждого из них. Пусть также задан подъём $P^l_1$ --- 
    $l$-мерный многогранник $\widetilde{P^l_1}$ в $\mathbb E^{d+1}$, тогда \emph{подъёмом $P^m_2$ до 
    $\widetilde{P^l_1}$} будем называть такой подъём $\widetilde{P^m_2}$ (в том же $\mathbb E^{d+1}$), 
    что пересечение $\widetilde{P^l_1}\bigcap\widetilde{P^m_2}$ является гранью каждого из 
    многогранников $\widetilde{P^l_1}$ и $\widetilde{P^m_2}$, и его проекция на $\mathbb E^d$ 
    совпадает с $P^l_1 \bigcap P^m_2$. Иными словами --- если $\widetilde{P^l_1}\bigcap
      \widetilde{P^m_2}$ является подъёмом $P^l_1 \bigcap P^m_2$.
  } \end{definition}
	
	\begin{lemma}\label{lemmaLiftFacet}~{
		Пусть заданы два выпуклых $d$-мерных многогранника $P_1, P_2 \subset \mathbb E^d$, смежные по 
		$(d-1)$-грани $F^{d-1}$ с нормалью $\mathbf{n} \in \mathbb E^d$ к этой грани. Пусть данное 
		пространство $\mathbb E^d$ вложено в $\mathbb E^{d+1}$ как гиперплоскость и задан некоторый 
		подъём многогранника $P_1$ --- $d$-мерный многогранник $\widetilde{P_1} \subset \mathbb E^{d+1}$ 
		с нормалью $\mathbf{n_1} \in \mathbb E^{d+1}$. Тогда существует единственный подъём многогранника
		$P_2$ до $\widetilde{P_1}$ с вектором нормали этого подъёма равным $\mathbf{n_1} + \mathbf{n}$.
	} \end{lemma}
	
	\smallskip
	{\it Доказательство.}
		Обозначим через $\Pr_1$ ограничение $\Pr|_{\aff{\widetilde{P_1}}}: \aff{\widetilde{P_1}} \to 
		\mathbb E^d$ проекции $\Pr$ на $d$-мерное аффинное подпространство многогранника 
		$\widetilde{P_1}$. Размерность образа и прообраза при проекции $\Pr_1$ совпадают, значит ядро 
		этого отображения нулевое и $\Pr_1$ биективно переводит $k$-мерные грани $\widetilde{P_1}$ в 
		$k$-мерные грани $P_1$ для произвольного $0 \leqslant k \leqslant d$. Кроме того, по 
		определению подъёма, $\mathbf{n}_1 \notparallel \mathbb E^d$. Обозначим через 
		$\widetilde{F^{d-1}}$ прообраз $F^{d-1}$ при проекции $\Pr_1$, а через $\widetilde{L^{d-1}}$ и 
		$L^{d-1}$ --- их аффинные оболочки. Очевидно $\Pr_1\parenth{\widetilde{L^{d-1}}} = L^{d-1}$.
		
		Покажем, что вектор $\mathbf{n_1} + \mathbf{n}$ ортогонален $\widetilde{L^{d-1}}$. Так как
		$\mathbf{n}_1 \bot \aff{\widetilde{P_1}}$, то $\mathbf{n}_1 \bot \widetilde{L^{d-1}} \subset
			\aff{\widetilde{P_1}}$. Пусть $\mathbf{n}_0 \in \mathbb E^{d+1}$ --- нормаль к $\mathbb E^d$, 
		тогда очевидно, что $\mathbf{n}_0$ отображается в нулевой вектор проекцией $\Pr$ всего 
		пространства $\mathbb E^{d+1}$ и является образующей в $\ker \Pr$. Отсюда и из того, что 
		$\Pr\parenth{\widetilde{L^{d-1}}} = L^{d-1}$, следует, что $\widetilde{L^{d-1}}$ принадлежит 
		$L^{d-1}\times\lin{\mathbf{n}_0}$, то есть $\widetilde{L^{d-1}}$ лежит в перпендикулярной к 
		$\mathbb E^d$ $d$-плоскости, содержащей $L^{d-1}$. По условию $\mathbf{n} \bot L^{d-1}$, кроме 
		того	$\mathbf{n} \bot \mathbf{n}_0$ так как $\mathbf{n}_0 \bot \mathbb E^d$ и $\mathbf{n} \in 
		  \mathbb E^d$ по условию. Следовательно $\mathbf{n} \bot L^{d-1}\times\lin{\mathbf{n}_0}$, 
		откуда имеем $\mathbf{n} \bot \widetilde{L^{d-1}}$. Таким образом, показано, что $\mathbf{n_1} 
		  + \mathbf{n}$ ортогонален $\widetilde{L^{d-1}}$.
		
		Выберем произвольную точку $A \in \widetilde{L^{d-1}}$ и проведём через неё $d$-мерную плоскость 
		$\pi^d_2$, ортогональную вектору $\mathbf{n_1} + \mathbf{n}$. По нормали и одной своей точке 
		гиперплоскость $\pi^d_2$ восстанавливается однозначно. По построению $\widetilde{L^{d-1}}$ 
		параллельно $\pi^d_2$ и их пересечение непусто, значит  $\widetilde{L^{d-1}}$ принадлежит 
		$\pi^d_2$. Отсюда имеем, что $\aff{\widetilde{P_1}} \bigcap \pi^d_2 = \widetilde{L^{d-1}}$, так
		как это пересечение содержит $\widetilde{L^{d-1}}$ и при этом не может быть $d$-мерным: нормали
		этих $d$-подпространств $\mathbf{n}_1$ и $\mathbf{n}_1 + \mathbf{n}$ не параллельны. Плоскость
		$\pi^d_2$ не ортогональна $\mathbb E^d$, так как $\parenth{\mathbf{n}_1 + \mathbf{n}} 
		  \notparallel E^d$. Значит ограничение проекции $\Pr$ на $\pi^d_2$ --- проекция 
		$\Pr_2 := \Pr|_{\pi^d_2}: \pi^d_2 \to \mathbb E^d$ также является невырожденным аффинным 
		преобразованием. В частности, для него однозначно определено обратное аффинное преобразование 
		$\Pr^{-1}_2 : \mathbb E^d \to \pi^d_2$.
		
		Убедимся, что $\Pr^{-1}_2\parenth{P_2}$ и есть искомый подъём $\widetilde{P_2}$. Отметим, что 
		ограничения $\Pr_1$ и $\Pr_2$ на $\widetilde{L^{d-1}}$ суть одно и то же преобразование --- 
		ограничение порождающей их проекции $\Pr$ на $\widetilde{L^{d-1}}$. Отсюда 
		$\Pr^{-1}_2\parenth{F^{d-1}} =  \widetilde{F^{d-1}}$. Далее, так как $\aff{\widetilde{P_1}} 
		  \bigcap \pi^d_2 = \widetilde{L^{d-1}}$, то пересечение $\Pr^{-1}_2\parenth{P_2}$ с 
		$\widetilde{P_1}$ принадлежит их общему подпространству $\widetilde{L^{d-1}}$. Тогда 
		$\Pr^{-1}_2\parenth{P_2} \bigcap \widetilde{P_1} = 
		  \Pr^{-1}_2\parenth{P_2} \bigcap \widetilde{L^{d-1}}\bigcap \widetilde{P_1} = 
		  \Pr^{-1}_2\parenth{P_2 \bigcap L^{d-1}} \bigcap \widetilde{P_1} = 
		  \Pr^{-1}_2\parenth{F^{d-1}} \bigcap \widetilde{P_1} = 
		  \widetilde{F^{d-1}} \bigcap \widetilde{P_1} = \widetilde{F^{d-1}}$.
		Таким образом, пересечение $\Pr^{-1}_2\parenth{P_2}$ с $\widetilde{P_1}$ является
		гранью каждого и является подъёмом для $F^{d-1} = P_1 \bigcap P_2$. По построению 
		$\Pr^{-1}_2\parenth{P_2}$ имеет размерность $d$ и ортогонален вектору $\mathbf{n}_1 + 
		  \mathbf{n}$. Значит $\Pr^{-1}_2\parenth{P_2}$ действительно является подъёмом $P_2$ до 
		$\widetilde{P_1}$. Кроме того, мы показали, что $d$-плоскость $\aff{\widetilde{P_2}}$ определена однозначно, так же как и 
		отображение $\Pr^{-1}|_{\aff{\widetilde{P_2}}}: \mathbb E^d \to \aff{\widetilde{P_2}}$. Значит 
		подъём $\widetilde{P_2}$ существует и определён однозначно.
	\hfill $\Box$ \par\bigskip
	
	Если, как в условии леммы \ref{lemmaLiftFacet}, подъём ячейки $P_2$ до поднятой ячейки 
	$\widetilde{P_1}$ задан при помощи вектора $\mathbf{n}$ ортогонального общей гиперграни $P_1$ и 
	$P_2$, то будем говорить, что $\widetilde{P_2}$ является подъёмом $P_2$ до $\widetilde{P_1}$
	\emph{с приращением нормали $\mathbf{n}$}.

\subsection{Построение женератрисы}\label{subsectionConstruction}

  Будем считать, что пространство $\mathbb E^d$, в котором задано разбиение $\mathcal T$, вложено 
  в пространство $\mathbb E^{d+1}$ как гиперплоскость $x_{d+1} = 0$. Зафиксируем нормаль 
  $\mathbf{n}_0 = (0, \ldots, 0, -1)$ к этой гиперплоскости. Для произвольных смежных 
  $d$-многогранников $P_i, P_j \subset \mathcal T$ с общей гипергранью $P_i\bigcap P_j$ через 
  $\mathbf{n}_{i, j}$ будем обозначать единичную нормаль к $P_i\bigcap P_j$, направленную от $P_i$ 
  к $P_j$, а через $s\parenth{P_i\bigcap P_j}$ --- значение заданной канонической нормировки $s$ на
  этой гиперграни разбиения. 
  
  Построим в $\mathbb E^{d+1}$ семейство $d$-многогранников $\mathcal G = \mathcal 
    G\parenth{\mathcal T, s}$ по следующим правилам:
  \begin{itemize} 
  	\item Каждый многогранник из $\mathcal G$ является подъёмом некоторой $d$-ячейки 
  	  $P \in \mathcal T$, занумерован при помощи некоторой цепи ячеек 
  	  $\brackets{P_0, P_{i_1}, \ldots, P_{i_k}, P}$, соединяющей $P_0$ c $P$, и обозначается 
  	  $L\brackets{P_0, P_{i_1}, \ldots, P_{i_k}, P}$
  	\item $L\brackets{P_0}$ совпадает с ячейкой $P_0$, нормаль к $L\brackets{P_0}$ (в $\mathbb 
  			E^{d+1}$) фиксирована и равна $\mathbf{n}_0$
  	\item $L\brackets{P_0, P_1, \ldots, P_{k-1}, P_k}$ является подъёмом $P_k$ до (ранее 
  		поднятого) $L\brackets{P_0, P_1, \ldots, P_{k-1}}$ с приращением нормали 
  		$s(P_{k-1}\bigcap P_k)\cdot\mathbf{n}_{k-1, k}$
  \end{itemize}
  
  Из леммы \ref{lemmaLiftFacet} следует, что для произвольной цепи $\brackets{P_0, \ldots, P_k}$
  каждый из $d$-многогранников $L\brackets{P_0, P_1}, L\brackets{P_0, P_1, P_2}, 
    L\brackets{P_0, \ldots, P_k}$ смежен с предыдущим (для $L\brackets{P_0, P_1}$ предыдущий ---
  $L\brackets{P_0}$ --- задан в определении) и по приращениям нормали восстанавливается однозначно.
  Будем называть многогранники $L\brackets{P_0, \ldots, P_k}$ гранями $\mathcal G$. Введём на
  $\mathcal G$ структуру полиэдрального комплекса:
  \begin{itemize}
  	\item Будем считать $d$-грани $L\brackets{P_0, P_{i_1}, \ldots, P_{i_m}}$ и 
  	  $L\brackets{P_0, P_{j_1},\ldots, P_{j_n}}$ совпадающими, если они совпадают геометрически как
  	  точечные множества в $\mathbb E^{d+1}$. В частности, отсюда вытекает необходимое условие
  	  $P_{i_m} = P_{j_n}$.
  	\item Две $d$-грани $\mathcal G$ считаем смежными тогда и только тогда, когда для них существует
  	  представление в виде $L\brackets{P_0, \ldots, P_{k-1}}$ и $L\brackets{P_0, \ldots, P_{k-1}, 
  	    P_k}$, где они следуют одна за другой в цепочке подъёмов.
  	\item Согласно определению \ref{definitionComplex}, заданные смежности $d$-мерных граней 
  	  $\mathcal G$ задают все остальные инцидентности в гранях меньших размерностей.
  \end{itemize}
  
  Далее под $\mathcal G$ мы будем понимать именно полиэдр с заданными инцидентностями (а не просто
  набор $d$-многогранников).

  \subsection{Однозначность проекции женератрисы}
  
  \emph{Женератрисой} разбиения $\mathcal T$ пространства $\mathbb E^d$ будем называть $d$-мерное 
  полиэдральное многообразие в объемлющем пространстве $\mathbb E^{d+1}$, ортогональная проекция
  которого на $\mathbb E^d = \aff{\mathcal T}$ однозначна и совпадает с разбиением $\mathcal T$ 
  (см. \cite{Voronoy, Davis, McMullen}).
  
  \begin{theorem}\label{theoremManifoldGeneratris}~ { 
    Полиэдр $\mathcal G_{\mathcal T}$, построенный по канонической нормировке $s$ является 
    женератрисой для $\mathcal T$. Обратно, если задана некоторая женератриса разбиения $\mathcal T$,
    то существует каноническая нормировка $s$ этого разбиения.    
  } \end{theorem}
  
  Доказательство этой теоремы разобьём на несколько лемм.

  \begin{lemma}\label{lemmaUniqueNormal}~{
    Все гиперграни полиэдра $\mathcal G_{\mathcal T}$, являющиеся подъёмами одной и той же 
    $d$-ячейки $P \subset \mathcal T$, параллельны друг другу.  
  } \end{lemma}
  
  \smallskip
  {\it Доказательство.}
    В определении $\mathcal G_{\mathcal T}$ однозначно задана (аддитивная) функция 
    приращения нормали при переходе от ячейки $P_{i-1}$ к смежной ячейке $P_i$:
    $g\brackets{P_{i-1}, P_i} = s\parenth{P_{i-1}\bigcap P_i}\mathbf{n}_{i-1, i}$. Воспользуемся 
    аддитивной версией теоремы \ref{theoremQualityTr}, чтобы проверить, что приращения задают 
    однозначную функцию $\mathop{\mathbf{n}}\parenth{P_i}$ нормали. Достаточно проверить, что 
    полное приращение вдоль любого $(d-2)$-примитивного пути $\brackets{P_{i_1}, P_{i_2}, 
    \ldots, P_{i_{t+1}} = P_{i_1}}$ равно нулю. То есть проверить условие 
    $\sum^t_{j=1}s\parenth{P_{i_j}\bigcap P_{i_{j+1}}}\mathbf{n}_{i_j,i_{j+1}} = \mathbf{0}$. 
    
    Данная цепь $(d-2)$-примитивна, значит все её ячейки содержат общую $(d-2)$-грань $F^{d-2}$. 
    Отметим, что приращение на участке цепи вида $\brackets{P_i,P_j,P_i}$ равно нулю. Отсюда 
    очевидно, что полное приращение равно приращению при обходе вокруг $F^{d-2}$, умноженному 
    на индекс цепи относительно такого обхода. Индекс мы определяем как количество раз, которое 
    была пройдена произвольная гипергрань $P_i\bigcap P_{i+1}$ данной цепи. При этом проходы в 
    направлении обхода вокруг $F^{d-2}$ считаем со знаком плюс, в обратном --- со знаком минус. 
    Приращение при обходе вокруг $F^{d-2}$ равно кручению нормировки $s$ вокруг $F^{d-2}$ и равно 
    нулю по определению канонической нормировки. Значит и полное приращение равно нулю. Таким 
    образом, функция $\mathbf{n}\parenth{P_i}$ определена однозначно, нормали всех подъёмов 
    $P_i$ в $\mathcal G$ равны, а сами подъёмы --- параллельны друг другу.
  \hfill $\Box$ \par\bigskip

  \begin{corollary}\label{corollaryCoincidenceCheck}~{ 
    Два подъёма $\widetilde{P}, \widetilde{P'} \subset \mathcal G_{\mathcal T}$ одного и того же 
    $d$-многогранника $P \subset \mathcal T$ представляют одну и ту же грань этого полиэдра тогда 
    и только тогда, когда имеют хотя бы одну общую точку.
  }\end{corollary}
  
  \begin{lemma}\label{lemmaCoherentProjection}~{ 
    Ортогональная проекция на $\mathcal T$ звезды произвольной $(d-2)$-грани $\widetilde{F^{d-2}}$ 
    полиэдра $\mathcal G_{\mathcal T}$ однозначна.
  } \end{lemma}

  \smallskip
	{\it Доказательство.}
		Предположим, проекция неоднозначна, тогда найдутся две различные $d$-грани $\widetilde{P}$ и 
		$\widetilde{P'}$ с пересекающимися по внутренней точке проекциями. Значит $\widetilde{P}$ и 
		$\widetilde{P'}$ являются подъёмами одной и той же $d$-ячейки разбиения $\mathcal T$. По условию
		эти два подъёма содержат общую грань $\widetilde{F^{d-2}}$. Согласно следствию 
		\ref{corollaryCoincidenceCheck} подъёмы $\widetilde{P}$ и $\widetilde{P'}$ совпадают, то есть 
		представляют одну и ту же $d$-грань полиэдра. Противоречие.
	\hfill $\Box$ \par\bigskip
	
	\smallskip
	{\it Доказательство теоремы \ref{theoremManifoldGeneratris}.}
	  Выберем внутри каждой ячейки $P_k \subset \mathcal T$ точку $p_k$. Пусть $P_{k-1}, P_k$ --- две
	  смежные ячейки в $\mathcal T$, $\widetilde{P_{k-1}}, \widetilde{P_k} \subset \mathcal 
	    G_{\mathcal T}$ --- два смежных подъёма $P_{k-1}, P_k$, а $\widetilde{p_{k-1}} \in 
	    \widetilde{P_{k-1}}, \widetilde{p_k} \in \widetilde{P_k}$ --- точки этих подъёмов, переходящие 
	  в $p_{k-1}$ и $p_k$ при ортогональной проекции. 
	  
	  Убедимся, что вектор $\widetilde{p_{k-1}p_k} = \widetilde{p_{k-1}}\widetilde{p_k}$ однозначно 
	  определён лишь выбором ячеек $P_{k-1}$ и $P_k$. Действительно, если $\overline{P_{k-1}}, 
	    \overline{P_k} \subset \mathcal G_{\mathcal T}$ --- два других их смежных подъёма, то по лемме
	  \ref{lemmaUniqueNormal}, имеем $\widetilde{P_{k-1}} \parallel \overline{P_{k-1}}$,
	  $\widetilde{P_k} \parallel \overline{P_k}$. Отсюда очевидно, что $\overline{P_{k-1}}\bigcup
	    \overline{P_k}$ является параллельным переносом $\widetilde{P_{k-1}}\bigcup \widetilde{P_k}$
	  вдоль нормали $\mathbf{n}_0 \bot \mathbb E^d$. Но тогда соответствующий вектор 
	  $\overline{p_{k-1}}\overline{p_k}$ равен своему параллельному переносу $\widetilde{p_{k-1}p_k}$.
	  
	  Таким образом определены векторные приращения $v\brackets{P_{k-1}, P_k} = 
	    \widetilde{p_{k-1}p_k}$ на парах смежных ячеек разбиения $\mathcal T$. Воспользуемся 
	  аддитивной версией теоремы \ref{theoremQualityTr}, чтобы показать, что такие приращения задают
	  однозначную вектор-функцию $v\parenth{P}$ на $d$-ячейках $P$ разбиения $\mathcal T$. Для этого 
	  достаточно показать, что на любом $(d-2)$-примитивном цикле в $\mathcal T$ приращение равно 
	  нулю. Рассмотрим $(d-2)$-грань $F^{d-2}$ и примитивный в ней цикл $\brackets{P_{i_1}, \ldots,
	    , P_{i_m} = P_{i_1}}$. Пусть $\widetilde{P_{i_1}}$ --- некоторая грань $\mathcal 
	    G_{\mathcal T}$ с проекцией $P_{i_1}$. Эта грань однозначно определяет подъём 
	  $\widetilde{F^{d-2}} \subset \widetilde{P_{i_1}}$ и соответствующую примитивную в 
	  $\widetilde{F^{d-2}}$ цепь $\brackets{\widetilde{P_{i_1}}, \ldots, \widetilde{P_{i_m}}}$. По
	  лемме \ref{lemmaCoherentProjection} ортогональная проекция 
	  $\mathop{St}_{\mathcal G_{\mathcal T}}\parenth{\widetilde{F^{d-2}}}$ на $\mathcal T$ однозначна.
	  Отсюда вытекает, что грани $\widetilde{P_{i_1}}$ и $\widetilde{P_{i_m}}$ совпадают, так же как и 
	  соответствующие точки $\widetilde{p_{i_1}}$ и $\widetilde{p_{i_m}}$. Приращение вдоль цепи
	  $\brackets{\widetilde{P_{i_1}}, \ldots, \widetilde{P_{i_m}}}$ равно 
	  $v\brackets{\widetilde{P_{i_1}}, \widetilde{P_{i_2}}} + v\brackets{\widetilde{P_{i_2}},
	    \widetilde{P_{i_3}}} + \ldots + v\brackets{\widetilde{P_{i_{m-1}}}, \widetilde{P_{i_m}}} = 
	    \widetilde{p_{i_1}p_{i_2}} + \ldots + \widetilde{p_{i_{m-1}}p_{i_m}} = 
	    \widetilde{p_{i_1}}\widetilde{p_{i_2}} + \ldots + \widetilde{p_{i_{m-1}}}\widetilde{p_{i_m}} =
	    \widetilde{p_{i_1}}\widetilde{p_{i_m}} = 0$
	    
	  Таким образом, приращения корректно задают некоторую вектор-функцию $v\parenth{P}$, определённую
	  однозначно, с точностью до выбора начального значения на одной фиксированной ячейке. Зададим 
	  значение $v\parenth{P_0} = p_0$. Тогда значение $v\parenth{P_k} = v\parenth{P_0} +
	    v\brackets{P_0, \ldots, P_k} = p_0 + \widetilde{p_0}\widetilde{p_1} + \ldots + 
	    \widetilde{p_{k-1}}\widetilde{p_k}$, где $\widetilde{p_i}$ --- это точка с проекцией $p_i$, 
	  лежащая в $d$-грани $\widetilde{P_i} \subset \mathcal G_{\mathcal T}$, поднятой вдоль данной 
	  цепи $\brackets{P_0, \ldots, P_k}$, а $\widetilde{p_0} = p_0$. Таким образом $v\parenth{P_k} =
	    p_0 + p_0\widetilde{p_k} = \widetilde{p_k}$. То есть значение $v\parenth{P_k}$, с одной 
	  стороны, определено однозначно, с другой --- принадлежит произвольному подъёму $\widetilde{P_k}
	  \subset \mathcal G_{\mathcal T}$  ячейки $P_k$. Согласно предложению 
	  \ref{corollaryCoincidenceCheck} все такие подъёмы совпадают и представляют одну и ту же грань 
	  $\mathcal G_{\mathcal T}$. Таким образом над каждой $d$-гранью $P$ из разбиения $\mathcal T$ 
	  есть в точности один её подъём, принадлежащий полиэдру $\mathcal G_{\mathcal T}$. Значит 
	  ортогональная проекция $\mathcal G_{\mathcal T}$ на $\mathbb E^d$ однозначна. В частности, 
	  отсюда следует, что $\mathcal G_{\mathcal T}$ является однолистной накрывающей $\mathbb E^d$, а 
	  значит является $d$-многообразием и женератрисой.
	  
	  Для доказательства теоремы в обратную сторону, достаточно выбрать для каждой $d$-грани 
	  женератрисы нормаль с $(d+1)$-ой координатой равной $-1$. Тем самым будут заданы 
	  приращения нормали (все будут параллельны $\mathbb E^d$) и нормировка $s$. Легко видеть, что 
	  она будет канонической.
	\hfill $\Box$ \par\bigskip

\subsection{Женератриса как функция}

  Из однозначности проекции $\mathcal G_{\mathcal T}$ на $\mathbb E^d$ следует, что 
  $\mathcal G_{\mathcal T}$ является графиком некоторой функции $G(x): \mathbb E^d \to \mathbb R$.
  Эту функцию, как и саму построенную поверхность, также будем называть женератрисой.

  \begin{lemma}\label{lemmaCompute}~{ 
    Пусть $P$ --- произвольная $d$-ячейка разбиения $\mathcal T$, $\brackets{P_0, \ldots, P_k = P}$
    --- некоторая цепь, соединяющая $P_0$ с $P$, тогда для любой пары точек $x_1, x_2 \in P$ 
    выполнено 
    $$G(x_1) - G(x_2) = (x_1 - x_2)^T\sum^{k-1}_{i=0}s\parenth{P_i\bigcap P_{i+1}}
      \mathbf{n}_{i,i+1}$$
  } \end{lemma}
  
  \smallskip
  \noindent {\it Доказательство.}
    По определению, точки женератрисы над ячейкой $P$ имеют вид $\parenth{x, G(x)}$, где $x \in P$.
    По построению, эти точки лежат в $d$-плоскости с нормалью $\mathbf{n}_0 + \sum^{k-1}_{i=0}
      s\parenth{P_i\bigcap P_{i+1}}\mathbf{n}_{i,i+1}$, где $\mathbf{n}_0 = \parenth{0, \ldots, 
      0, -1}^T \in \mathbb E^{d+1}$, $\mathbf{n}_{i,i+1} \parallel \mathbb E^d$. Из ортогональности имеем:
      $$0 = \parenth{\parenth{x_1, G(x_1)} - \parenth{x_2, G(x_2)}}^T\cdot
       \parenth{\sum^{k-1}_{i=0}s\parenth{P_i\bigcap P_{i+1}}\mathbf{n}_{i,i+1} + \mathbf{n}_0} =$$
      $$ = (x_1 - x_2)^T\cdot\sum^{k-1}_{i=0}s\parenth{P_i\bigcap P_{i+1}}\mathbf{n}_{i,i+1} +
        \parenth{G(x_1) - G(x_2)}\cdot(-1)$$
      Что и требовалось доказать.
  \hfill $\Box$ \par\bigskip
  
  \begin{corollary}\label{corollaryComputeEasy}~{ 
    Пусть задана точка $x \in P \subset \mathcal T$. Пусть также $\brackets{P_0, P_1, \ldots, 
      P_k = P}$ --- произвольная цепь, соединяющая начальную ячейку построения женератрисы $P_0$ с
      $P$, а точки $x_1, \ldots, x_k$ выбраны так, что $x_i \in P_{i-1}\bigcap P_i$. Тогда 
      $$G(x) = x^T\sum^{k-1}_{i=0}{s\parenth{P_i\bigcap P_{i+1}}\mathbf{n}_{i,i+1}} - 
        \sum^{k-1}_{i=0}{x^T_{i+1} s\parenth{P_i\bigcap P_{i+1}}\mathbf{n}_{i,i+1}}$$
  } \end{corollary}
  
  \smallskip
  \noindent {\it Доказательство.}
    По построению $P_0$ лежит в плоскости $x^{d+1} = 0$ в расширенном пространстве 
    $\mathbb E^{d+1}$. Значит $G(x_1) = 0$. Тогда $G(x) = G(x) - G(x_1) = \parenth{G(x) - G(x_k)}
      + \sum^{k-1}_{i=1}\parenth{G(x_{i+1}) - G(x_i)}$. Подставим выражения для этих разностей из
    леммы \ref{lemmaCompute} и приведём подобные.
  \hfill $\Box$ \par\bigskip
  
  Будем пользоваться стандартным обозначением $\epi G$ для надграфика женератрисы $\epi G = 
    \braces{(x, y) \in \mathbb E^d \times \mathbb R|~y \geqslant G(x)} \subset \mathbb E^{d+1}$.

	\begin{lemma}\label{lemmaConvexRidges}~{ 
	  Двугранный угол при произвольной $(d-1)$-грани надграфика $\epi G$ меньше развёрнутого.
	} \end{lemma}
	
	\smallskip
  \noindent {\it Доказательство.}
    Рассмотрим произвольный двугранный угол надграфика, образованный смежными подъёмами 
    $\widetilde{P}, \widetilde{P'} \subset \mathcal G_{\mathcal T}$ $d$-мерных ячеек $P, P' \subset
      \mathcal T$. При построении женератрисы для $\widetilde{P}$ была задана нормаль 
    $\mathbf{h}_P$ представленная в каноническом виде $\mathbf{h}_P = \mathbf{n}_0 + 
      \sum s\parenth{P_i\bigcap P_{i,i+1}}\mathbf{n}_{i,i+1}$, где сумма конечна и вычисляется 
    вдоль некоторой цепи в $\mathcal T \subset \mathbb E^d$, $\mathbf{n}_{i,i+1} \parallel
      \mathbb E^d$, а $\mathbf{n}_0 = \parenth{0,\ldots, 0, -1}^T$ --- вектор с единственной $-1$
    в $(d+1)$-ом разряде. Таким образом, $(d+1)$-ая координата $\mathbf{h}_P$ также равна $-1$, и 
    для $\widetilde{P}$ можно корректно определить ``верхнее'' и ``нижнее'' полупространства 
    относительно $\aff{\widetilde{P}}$: ``нижнее'' --- то, в которое выходит $\mathbf{h}_P$.
    У $\widetilde{P}$ и $\widetilde{P'}$ общая $(d-1)$-грань, поэтому достаточно показать, что 
    хотя бы одна точка $\widetilde{P'}$ лежит в ``верхнем'' полупространстве (без границы) 
    относительно $\aff{\widetilde{P}}$. Из этого будет следовать, что вся грань $\widetilde{P'}$
    (кроме точек $\widetilde{P}\bigcap\widetilde{P'}$) лежит в этом полупространстве. Это и будет
    означать, что двугранный угол между ними (относящийся к $\epi G$), меньше развёрнутого.
    
    Рассмотрим произвольный отрезок $x_0x_2 \subset \mathbb E^d$ такой, что $x_0 \in \int{P}, 
      x_2 \in \int{P'}$ и отрезок пересекает гипергрань $P\bigcap P'$ в точке $x_1$. По лемме
    \ref{lemmaCompute} имеем $G(x_1) - G(x_0) = (x_1 - x_0)^T(\mathbf{h}_P - \mathbf{n}_0)$.
    Очевидно, что уравнение $y - G(x_0) = (x - x_0)^T(\mathbf{h}_P - \mathbf{n}_0)$ задаёт 
    $d$-плоскость $\aff{\widetilde{P}}$. Значит для точки $(x_2, y_2) \in \aff{\widetilde{P}}$
    имеем $y_2 - G(x_0) = (x_2 - x_0)^T(\mathbf{h}_P - \mathbf{n}_0)$. Из той же леммы 
    \ref{lemmaCompute} имеем $G(x_2) - G(x_1) = (x_2 - x_1)^T(\mathbf{h}_P - \mathbf{n}_0 + 
      s\parenth{P\bigcap P'}\mathbf{n}')$, где $\mathbf{n}'$ --- единичная нормаль к $P\bigcap P'$,
    направленная от $P$ к $P'$. Отсюда $G(x_2) = y_2 + (x_2 - x_1)^Ts\parenth{P\bigcap P'}
      \mathbf{n}'$. Вектор $(x_2 - x_1)$ пересекает $P\bigcap P'$ и направлен строго внутрь $P'$, 
    Нормаль $\mathbf{n}'$ к $P\bigcap P'$ также направлена внутрь $P'$, $s\parenth{P\bigcap P'} 
      > 0$. Значит $(x_2 - x_1)^Ts\parenth{P\bigcap P'}\mathbf{n}' > 0$ и $G(x_2) > y_2$. 
    Лемма доказана. 
  \hfill $\Box$ \par\bigskip

  \begin{theorem}\label{theoremGeneratrissConvex}~{ 
    Женератриса $\mathcal G_{\mathcal T}$ разбиения $\mathcal T$, построенная по канонической 
    нормировке $s$, задаёт выпуклую непрерывную кусочно-линейную функцию $G(x)$ и, соответственно, 
    является границей выпуклого $(d+1)$-мерного полиэдра.
  } \end{theorem}

  \smallskip
  \noindent {\it Доказательство.} 
    Кусочная линейность $G(x)$ следует из определения женератрисы, непрерывность --- из теоремы 
    \ref{theoremManifoldGeneratris}. Остаётся показать выпуклость $G(x)$, то есть что для любой 
    пары точек $x, y \in \mathbb E^d$ отрезок, соединяющий точки $\parenth{x, G(x)}$ и 
    $\parenth{y, G(y)}$, принадлежит надграфику $\epi{G(x)}$. Для любого малого $\varepsilon > 0$ в 
    $\varepsilon$-окрестности $x$ и $y$ найдутся, 
    соответственно, точки $x'$ и $y'$, которые лежат строго внутри некоторых $d$-ячеек разбиения 
    $\mathcal T$, прямая $x'y'$ не пересекает $(d-2)$-остов разбиения $\mathcal T$, а $(d-1)$-грани
    пересекает не более чем по одной точке. Действительно, выберем сначала в $\frac{\varepsilon}{4}$-
    окрестностях $x$ и $y$ точки $x'''$ и $y'''$, не принадлежащие $\mathcal F^{d-1}$. Это можно 
    сделать, так как $(d-1)$-остов является локально-конечным объединением множеств $d$-мерной 
    меры нуль в $\mathbb E^d$. В $\frac{\varepsilon}{4}$-окрестностях $x'''$ и $y'''$ выберем точки 
    $x''$ и $y''$ так, что направление $x''y''$ не параллельно ни одной $(d-1)$-грани разбиения 
    $\mathcal T$. Это возможно, так как таких гиперграней счётное количество, каждая из них 
    ``запрещает'' направления некоторого проективного $(d-2)$-подпространства $P\mathbb R^{d-1}$. 
    По теореме Бэра, объединение счётного числа таких подмножеств меры нуль не покрывает никакое 
    открытое множество в проективном пространстве $P\mathbb R^d$ всех направлений в $\mathbb E^d$, 
    в том числе --- никакую окрестность направления $x'''y'''$. Наконец, спроектируем все грани 
    $(d-2)$-остова $\mathcal K^{d-2}$ вдоль направления $x''y''$ в 
    дополнительное $(d-1)$-подпространство (относительно $\mathbb E^d$). Грани остова перейдут в 
    счётное число многогранников размерности не более $(d-2)$, то есть множеств меры нуль в этом 
    $(d-1)$-подпространстве. Это множество не покрывает $\frac{\varepsilon}{4}$-окрестности точки, в 
    которую проектируется сама прямая $x''y''$. Значит, найдётся параллельный перенос прямой на 
    вектор длины не более $\frac{\varepsilon}{4}$, который не пересекает $\mathcal K^{d-2}$. Образы 
    точек $x''$ и $y''$ при этом переносе обозначим $x'$ и $y'$, которые и есть искомые точки.
    
    Рассмотрим двумерную плоскость $\pi \in \mathbb E^{d+1}$, которая содержит прямую $x'y'$ и
    ортогональна $\mathbb E^d$. Для любой точки $z$ прямой $x'y'$, соответствующая точка 
    женератрисы $\parenth{z, G(z)}$ также принадлежит $\pi$. По построению $\pi$ пересекает только 
    $d$ и $(d-1)$-мерные грани $\mathcal T$, причём $(d-1)$-мерные --- только в одной точке. 
    Отсюда следует, что то же самое верно и для пересечения $\pi$ с $\mathcal G_{\mathcal T}$.
    
    Множество $M$ называется \emph{локально выпуклым} в точке $x$, если найдётся такое открытое 
    множество $U\ni x$, что $M \bigcap U$ выпукло. Покажем, что надграфик $\epi{G|_{x'y'}}$ 
    функции $G(x)$, ограниченной на прямую $x'y'$, является локально выпуклым. Для внутренних точек
    надграфика и внутренних точек одномерных рёбер надграфика это очевидно. Вершины 
    $\epi{G|_{x'y'}}$ --- пересечение $\pi$ с внутренностью $(d-1)$-мерных граней $\mathcal 
      G_{\mathcal T}$. Согласно лемме \ref{lemmaConvexRidges}, двугранный угол в $(d-1)$-гранях
    $\mathcal G_{\mathcal T}$ меньше развёрнутого. Отсюда следует, что надграфик $\epi{G}$
    локально выпуклый во внутренних точках $(d-1)$-граней: достаточно рассмотреть шаровую 
    окрестность соответствующей точки, не пересекающую других граней $\mathcal G_{\mathcal T}$.
    Таким образом $\epi{G|_{x'y'}}$ в любой своей вершине является пересечением локально выпуклого
    в этой точке $\epi{G}$ и выпуклой двумерной плоскости. Отсюда $\epi{G|_{x'y'}}$
    также является локально выпуклым в вершинах, а значит и во всех своих точках.
    
    Воспользуемся теоремой Бёрдона \cite[Теорема 7.5.1]{Beardon}  о том, что замкнутое 
    линейно-связное локально выпуклое подмножество двумерной плоскости является выпуклым. Из 
    теоремы следует, что $\epi{G|_{x'y'}}$ выпуклый. Значит ему принадлежит весь отрезок 
    соединяющий $\parenth{x', G(x')}$ и $\parenth{y', G(y')}$. Значит этот отрезок принадлежит и 
    надграфику $\epi{G}$. Так как $G$ непрерывна, а точки $x', y'$ выбирались сколь угодно 
    близко к $x, y$, то и отрезок, соединяющий $\parenth{x, G(x)}$ и $\parenth{y, G(y)}$, также 
    принадлежит $\epi{G}$. Теорема доказана.
  \hfill $\Box$ \par\bigskip

  \begin{corollary}\label{corollaryNonnegative}~{
    Функция женератрисы $G(x)$ неотрицательна на $\mathbb E^d$.
  } \end{corollary}
  
  \smallskip
  \noindent {Доказательство.} 
    По доказанному, надграфик $\epi{G(x)}$ выпуклый. Гипергранью этого надграфика является 
    $d$-мерная ячейка $P_0$ --- начальная ячейка построения женератрисы. Из выпуклости следует, что
    $\epi{G(x)}$ лежит в одном полупространстве относительно гиперплоскости $\aff{P_0} = 
      \braces{x^{d+1} = 0}$.
  \hfill $\Box$ \par\bigskip

\section{Женератриса Вороного}
\subsection{Вычисление значений}

  Вернёмся к разбиениям на параллелоэдры. Пусть в пространстве $\mathbb E^d$ задано нормальное
  разбиение $\mathcal T_P$, порождённое параллелоэдром $P$, и существует каноническая нормировка $s$
  этого разбиения. По определению параллелоэдра, каждая ячейка этого разбиения совмещается с каждой
  другой при помощи некоторого параллельного переноса. Легко установить, что такой параллельный 
  переводит разбиение $\mathcal T_P$ в себя \cite{Dolbilin_new}. Отсюда, как уже упоминалось, 
  следует, что центры параллелоэдров разбиения образуют $d$-мерную целочисленную решётку 
  $\Lambda^d(P)$ для некоторого базиса в $\mathbb E^d$. Вектора этой решётки --- в точности вектора 
  всех параллельных переносов $\mathbb E^d$, сохраняющих разбиение $\mathcal T_P$. Если при всех таких
  переносах нормировка $s$ на гипергранях $\mathcal T_P$ сохраняется, то будем называть её 
  \emph{трансляционно инвариантной}. Доказательство следующей леммы было предложено А. Гарбером в 
  личных обсуждениях.
  
  \begin{lemma}\label{lemmaInvariantScaling}~{
    Если для данного разбиения $\mathcal T_P$ существует некоторая каноническая нормировка $s$, то 
    существует и трансляционно инвариантная каноническая нормировка данного разбиения.
  } \end{lemma}
  
  \smallskip
	{\it Доказательство.}
	  Выберем произвольный параллелоэдр $P_0$ данного разбиения. Если нормировки на какой-то его паре
	  противоположных гиперграней $F_1$ и $F'_1$ не равны, то из леммы \ref{lemmaCross} следует, что
	  все $(d-2)$-подграни $F_1$ стандартны: если есть хотя бы одна примитивная, то ей соответствует 
	  6-поясок, содержащий $F_1$ и $F'_1$ и согласно этой лемме, $s(F_1) = s(F'_1)$, что неверно. 
	  Тогда заменим нормировки всех гиперграней в $\mathcal T_P$ параллельных $F_1$ на $s(F_1)$. 
	  Легко видеть, что все эти грани --- суть параллельные переносы $F_1$, при которых разбиение 
	  переходит в себя. Значит, по доказанному, такая замена нормировки затронет лишь звёзды 
	  стандартных $(d-2)$-граней. Из предложения \ref{propositionDminus2Scaling} следует, что в каждой
	  из этих $(d-2)$-граней новая нормировка осталась канонической (то есть кручение вокруг данной 
	  $(d-2)$-грани равно 0). Значит и в целом новая нормировка --- каноническая.
	  
	  Таким образом, можем считать, что нормировка $s$ совпадает на противоположных гранях $P_0$. 
	  Значит параллельные переносы $P_0$ вместе с его нормировкой на вектора решётки $\Lambda^d(P)$
	  корректно задают некоторую трансляционно инвариантную нормировку $s'$. Покажем, что $s'$ также 
	  каноническая. Кручение $\Delta_{s'}\parenth{F^{d-2}}$ вокруг произвольной $(d-2)$-грани 
	  разбиения равно кручению нормировки $s'$ вокруг некоторого $\parenth{F^{d-2}}'$ --- 
	  параллельного переноса $F^{d-2}$, который принадлежит $P_0$. Поэтому сразу считаем, что $F^{d-2} 
	    \subset P_0$. Если $F^{d-2}$ стандартна, то
	  из предложения \ref{propositionDminus2Scaling} следует $\Delta_{s'}\parenth{F^{d-2}} = 0$. Если
	  $F^{d-2}$ примитивная, то обозначим через $F_2$ и $F_3$ гиперграни $P_0$, содержащие $F^{d-2}$, 
	  $\brackets{F_1, F_2, \ldots, F_6}$ --- соответствующий 6-поясок, $F_0$-третья гипергрань 
	  разбиения, содержащая $F^{d-2}$. По лемме \ref{lemmaCross} для исходной нормировки $s$ и 
	  накрест лежащих гиперграней $F_0$ и $F_1$ имеем $s(F_1) = s(F_0)$. По определению $s'$ 
	  выполнено $s'(F_0) = s(F_1) = s(F_0),~s'(F_2) = s(F_2),~s'(F_3) = s(F_3)$. Кручение нормировки 
	  $s$ вокруг $F^{d-2}$ равно 0, и его слагаемые в точности равны слагаемым кручения нормировки 
	  $s'$ вокруг $F^{d-2}$. Значит $\Delta_{s'}\parenth{F^{d-2}} = \Delta_{s}\parenth{F^{d-2}} = 0$. 
	  Значит $s'$ --- каноническая и трансляционно инвариантная. 
	\hfill $\Box$ \par\bigskip
	
	Далее считаем, что заданная на разбиении $\mathcal T_P$ каноническая нормировка $s$ трансляционно 
	инвариантна. Согласно приведённой в разделе \ref{subsectionConstruction} конструкции, 
	$\mathcal T_P$ и $s$ однозначно задают 
	женератрису $\mathcal G_P \subset \mathbb E^{d+1}$. \emph{Фасетным вектором} разбиения $\mathcal 
	  T_P$ называется всякий вектор, соединяющий центры двух смежных по общей гиперграни ячеек. 
	Зафиксируем некоторую ячейку разбиения $P_0$. Пусть $\braces{P_1, P_2, \ldots P_k}$ --- все 
	параллелоэдры, смежные с $P_0$. Центр ячейки $P_i$ обозначим $c(P_i)$. Очевидно, все фасетные 
	вектора разбиения можно представить в виде $\mathbf{p}_i = c(P_0)c(P_i)$ (все эти вектора делятся 
	на пары дающих в сумме ноль). Также обозначим $\mathbf{m}_i = s(P_0 \bigcap P_i)\mathbf{n}_{0,i}$. 
	Через $P(\p)$ будем обозначать копию параллелоэдра $P$ с центром в точке $\p$.

	\begin{lemma}\label{lemmaSymmetry}~{ 
	  Для произвольных индексов $i, j$ выполнено равенство $\p^T_i\m_j = \p^T_j\m_i$.
	} \end{lemma}

  \smallskip
  \noindent {\it Доказательство.} 
    Вычислим значение $G(\p_1 + \p_2)$ двумя способами: вдоль цепи $\brackets{P_0, P(\p_1), P(\p_1 
      + \p_2)}$ и вдоль цепи $\brackets{P_0, P(\p_2), P(\p_1 + \p_2)}$. Середины отрезков
    $(0, \p_1)$ и $(\p_1, \p_1 + \p_2)$ являются центрами симметрии соответствующих пар смежных 
    параллелоэдров и принадлежат соответствующим гиперграням. Выберем эти точки $x_1 = 
      \frac{\p_1}{2}$ и $x_2 = \p_1 + \frac{\p_2}{2}$ как вспомогательные, положим $x = \p_1 + 
      \p_2$ и воспользуемся следствием \ref{corollaryComputeEasy}. Получаем 
    $$G(\p_1 + \p_2) = (\p_1 + \p_2)^T(\m_1 + \m_2) - \parenth{\frac{\p_1}{2}^T\m_1 + \parenth{\p_1
      + \frac{\p_2}{2}}^T\m_2}$$
    Тут мы воспользовались наблюдением, что гипергрань между $P(\p_1)$ и $P(\p_1 + \p_2)$ параллельна
    гиперграни между $P(0)$ и $P(\p_2)$, значит их нормали параллельны, а нормировки равны. 
    Аналогично для второй цепи имеем 
    $G(\p_1 + \p_2) = (\p_1 + \p_2)^T(\m_1 + \m_2) - \parenth{\frac{\p_2}{2}^T\m_2 + \parenth{\p_2
      + \frac{\p_1}{2}}^T\m_1}$.
    Приводим подобные и получаем требуемое равенство.
  \hfill $\Box$ \par\bigskip
  
  Пусть $\mcP = \parenth{\p_1 | \ldots | \p_k}$ --- это матрица $d\times k$, составленная из 
  всех фасетных векторов параллелоэдра $P_0$, где $i$-й столбец представлен вектором $\p_i$. Матрица
  $\mcM = \parenth{\m_1 | \ldots | \m_k}$ --- матрица, составленная по тем же правилам из 
  векторов $\m_i$.
  
  \begin{lemma}\label{lemmaGPFormula}~{ 
    Для произвольных целых чисел $l_1, \ldots, l_k$ и вектора $L = \parenth{l_1, \ldots, l_k}$
    выполнено:
    $$G\parenth{l_1\p_1 + \ldots + l_k\p_k} = \frac12\parenth{l_1\p_1 + \ldots + l_k\p_k}^T
      \parenth{l_1\m_1 + \ldots + l_k\m_k} = \frac12L\mcP^T\mcM L^T$$
  } \end{lemma}
  
  \smallskip
  \noindent {\it Доказательство.} 
    Рассмотрим вспомогательные точки $x_1 = \frac{\p_1}{2}, x_2 = \p_1 + \frac{\p_1}{2}, \ldots,
      x_{l_1} = (l_1 - 1)\p_1 + \frac{\p_1}{2}, x_{l_1 + 1} = l_1\p_1 + \frac{\p_2}{2}, \ldots, 
      x_{\sum{l_i}} = l_1\p_1 + \ldots + (l_k - 1)\p_k + \frac{\p_k}{2}$. Если некоторое $l_i$ 
    отрицательно, то рассматриваем слагаемое $l_i\p_i$ как $|l_i|(-\p_i)$. Воспользуемся 
    следствием \ref{corollaryComputeEasy}. Первое равенство из доказываемой цепочки получается 
    приведением подобных и использованием леммы \ref{lemmaSymmetry}. Второе равенство --- 
    матричная запись полученного выражения.
  \hfill $\Box$ \par\bigskip
  
  \begin{remark} {
    Из следствия \ref{corollaryComputeEasy} напрямую вытекает, что вектор
    $\parenth{l_1\m_1 + \ldots + l_k\m_k} = \mcM L^T$ является градиентом 
    $G$ над ячейкой $P\parenth{l_1\p_1 + \ldots + l_k\p_k}$.
  } \end{remark}

  Элегантное доказательство следующей леммы практически дословно повторяет доказательство из 
  работы \cite{DezaGrishukhin}. Для полноты изложения и самодостаточности работы мы приводим его 
  в терминах канонических нормировок.
  
  \begin{lemma}\label{lemmaMatrixEquivalence}~{
    Существует единственная симметричная невырожденная $d\times d$ матрица $Q$ такая, что 
    $\mcM = Q\mcP$.
  } \end{lemma}

  \smallskip
  \noindent {\it Доказательство.}
    Выберем набор индексов $\braces{i_1, \ldots, i_d} \subset \braces{1, \ldots, k}$ такой, что 
    вектора $\p_{i_1}, \ldots, \p_{i_d}$ линейно независимы. Такой набор из $d$ фасетных векторов 
    обязательно найдётся, так как, по очевидным причинам, полный набор $\p_1, \ldots, \p_k$ своими 
    целочисленными линейными комбинациями порождает решётку $\Lambda^d(P)$ \cite{Dolbilin_new}. 
    Обозначим через $\mcP_0$ минор $\mcP$ составленный только из столбцов $\p_{i_1}, \ldots, 
      \p_{i_d}$, $\mcM_0$ --- минор $\mcM$, составленный из столбцов с теми же индексами. По выбору
    $\mcP_0$ --- невырожденная матрица.
      
    Из равенства $\p^T_i \m_j = \p^T_j \m_i$ следует, что $\mcP^T\mcM = \mcM^T\mcP$ и $\mcP^T_0\mcM 
      = \mcM^T_0\mcP$. Обозначим $Q := \parenth{\mcP^T_0}^{-1}\mcM^T_0$, тогда $\mcM = 
      \parenth{\mcP^T_0}^{-1}\mcM^T_0\mcP = Q\mcP$. Отсюда, в частности, $\m_i = Q\p_i$ и 
    следовательно $\mcM_0 = Q\mcP_0$. Значит для $Q$ также выполнено $Q = \mcM_0\mcP^{-1}_0$, 
    откуда $Q = \parenth{\mcP^T_0}^{-1}\mcM^T_0 = \parenth{\mcM_0\mcP^{-1}_0}^T = Q^T$, то есть $Q$ 
    симметрична.
    
    Из равенства $\mcM = Q\mcP$ и из того, что ранг матрицы $\mcM$ (нормалей к гиперграням 
    параллелоэдра $P$) равен $d$, следует, что ранг $d\times d$ матрицы $Q$ также равен $d$ и $Q$
    невырождена. Единственность следует из необходимого равенства $\mcM_0 = Q\mcP_0$.
  \hfill $\Box$ \par\bigskip

  \begin{corollary}\label{corollaryValueOnCenters}~{
    Пусть $x$ --- центр произвольной ячейки $P(x)$ разбиения $\mcT_P$. Тогда $G(x) = \frac12x^tQx$.
  } \end{corollary}

  \smallskip
  \noindent {\it Доказательство.}
    Центр произвольного параллелоэдра разбиения представляется в виде некоторой целочисленной 
    линейной комбинации $x = l_1\p_1 + \ldots + l_k\p_k = \mcP L^T$. По лемме \ref{lemmaGPFormula} 
    $G(x) = G\parenth{l_1\p_1 + \ldots + l_k\p_k} = \frac12L\mcP^T\mcM L^T$, что по лемме
    \ref{lemmaMatrixEquivalence} равно $\frac12\parenth{L\mcP^T}Q\mcP L^T = \frac12x^TQx$.
  \hfill $\Box$ \par\bigskip
  
  Обозначим через $Q(x)$ квадратичную форму на $\mathbb E^d$, заданную симметричной матрицей $Q$:
  $Q(x) = \frac12x^tQx$. Непосредственно из определения и следствия \ref{corollaryValueOnCenters}
  следует
  
  \begin{proposition}\label{propositionCoinside}~{
    Значения $Q(x)$ и $G(x)$ на центрах параллелоэдров разбиения $\mcT_P$ совпадают.
  } \end{proposition}
  
  \begin{lemma}\label{lemmaPositiveDefinite}~{
    Форма $Q(x)$ положительно определена.
  } \end{lemma}
  
  \smallskip
  \noindent {\it Доказательство.}
    Докажем сначала, что $Q(x)$ неотрицательно определена. Действительно, если для некоторой точки
    $x^*$ выполнено $Q(x^*) < 0$, то существует не содержащая нуля шаровая окрестность 
    $B_{\delta}(x^*)$, на точках которой $Q(x)$ также отрицательна. Значит $Q(x)$ отрицательна на 
    всех точках конуса над $B_{\delta}(x^*)$ с вершиной в нуле. Этот конус содержит шар сколь 
    угодно большого радиуса, а значит содержит некоторую точку $x_{\Lambda}$ решётки 
    $\Lambda^d(P)$. Тогда, с одной стороны, $Q(x_{\Lambda}) < 0$, с другой $Q(x_{\Lambda}) = 
      G(x_{\Lambda}) \geqslant 0$ по следствию \ref{corollaryNonnegative}. Противоречие. Значит
    $Q(x)$ неотрицательна определена. Кроме того, по лемме \ref{lemmaMatrixEquivalence}, матрица
    $Q$ невырождена, значит $Q(x)$ положительно определена.
  \hfill $\Box$ \par\bigskip
  
  \begin{theorem}\label{theoremTouch}~{
    График квадратичной формы $Q(x)$ --- эллиптический параболоид $\Pi_P \subset \mathbb E^{d+1}$.
    $\Pi_P$ вписан в женератрису $\mcG_P$ и касается её над всеми центрами ячеек разбиения $\mcT_P$.
  } \end{theorem}
  
  \smallskip
  \noindent {\it Доказательство.}
    То, что $\Pi_P$ является эллиптическим параболоидом, эквивалентно лемме 
    \ref{lemmaPositiveDefinite}. Значения $Q(x)$ и $G(x)$ совпадают над центрами ячеек $\mcT_P$.
    Докажем касание их графиков в этих точках. Для этого достаточно доказать пропорциональность 
    градиентов в центре $x$ произвольной ячейки. Пусть $x$ выражается через фасетные вектора 
    в виде целочисленной линейной комбинации $x = l_1\p_1 + \ldots + l_k\p_k$. Согласно замечанию
    к лемме \ref{lemmaGPFormula}, градиент $G(x)$ над ячейкой $P(x)$ равен $l_1\m_1 + \ldots + 
      l_k\m_k$. Из леммы \ref{lemmaMatrixEquivalence} это выражение равно $l_1Q\p_1 + \ldots + 
      l_kQ\p_k = Qx$. Градиент $Q(x)$ равен $Qx$ во всех точках. Касание доказано. Осталось 
    заметить, что эллиптический параболоид --- строго выпуклое $(d+1)$-мерное тело и лежит строго
    в одном полупространстве относительно любой своей касательной плоскости.
  \hfill $\Box$ \par\bigskip

\subsection{Искомое аффинное преобразование}

  По теореме из линейной алгебры, существует невырожденное линейное преобразование $\mcA_d$, при 
  котором положительно определённая квадратичная форма $Q$ переходит в стандартную квадратичную 
  форму с единичной матрицей $I = \mathop{diag}\parenth{1, 1, \ldots, 1}$. Обычно это утверждение
  формулируется в терминах замен координат и подразумевает, что квадратичной форме $Q(x)$ 
  соответствует такая форма $Q_{\mcA}$, что для всех $x \in \mathbb E^d$ верно 
  $Q_{\mcA}\parenth{\mcA_d x} = Q(x)$. Мы покажем, что линейное преобразование $\mcA_d$, с 
  точностью до дополнительных параллельных переносов, является искомым аффинным преобразованием,
  переводящим исходный параллелоэдр $P$ в некоторый параллелоэдр Вороного.
  
  Отображение $\mcA_d: \mathbb E^d \to \mathbb E^d$ продолжается до невырожденного линейного 
  преобразования $\mcA_{d+1}: \mathbb E^{d+1} \to \mathbb E^{d+1}$ по правилу 
  $\mcA_{d+1}: \column{x}{x^{d+1}} \mapsto \column{\mcA_d x}{x^{d+1}}$, для 
  произвольной $x \in \mathbb E^d$ и $x^{d+1} \in \mathbb R$. Тогда $\mcA_d = 
    \mcA_{d+1}|_{x^{d+1} = 0}$.
  
  Отображение $\mcA_d$ переводит разбиение $\mcT_P$ пространства $\mathbb E^d$ в разбиение
  $\mcA_d\mcT_P$, ячейки которого совмещаются каждая с каждой при помощи подходящего параллельного
  переноса. Значит многогранник $\mcA_dP$ также является параллелоэдром и 
  $\mcA_d\mcT_P = \mcT_{\mcA_dP}$.
  
  \begin{lemma}\label{lemmaMappingCommutes}~{
    Отображение $\mcA_{d+1}$ переводит параболоид $\Pi_P$ в стандартный параболоид $\Pi_I = 
      \braces{y = \parenth{y^1, \ldots, y^{d+1}}^T|~y^{d+1} = \sum^d_{i=1}\parenth{y^i}^2}$, 
    женератрису $\mcG_P$ --- в женератрису $\mcG_{\mcA_dP}$ разбиения $\mcT_{\mcA_dP}$. При этом
    $\Pi_I$ вписан в $\mcG_{\mcA_dP}$, а проекции точек касания --- суть центры параллелоэдров
    разбиения $\mcT_{\mcA_dP}$.
  } \end{lemma}
  
  \smallskip
  \noindent {\it Доказательство.}
    По определению $\mcA_{d+1}\column{x}{Q(x)} = \column{\mcA_dx}{Q(x)} = \column{\mcA_dx}
      {Q_{\mcA}\parenth{\mcA_dx}}$. Значит $\mcA_{d+1}\Pi_P \subset \Pi_I = \braces{
      \column{y}{Q_{\mcA}(y)}|~y\in\mathbb E^d}$. Так как $\mcA_d$ невырождено, то для 
    произвольного $y$ найдётся прообраз $x$. Значит $\mcA_{d+1}\Pi_P = \Pi_I$. Из определения
    следует, что $\mcA_{d+1}$ переводит проекции фигур на $\mathbb E^d$ в проекции образов этих
    фигур на $\mathbb E^d$. Отсюда следуют остальные утверждения леммы.
  \hfill $\Box$ \par\bigskip
  
  \begin{lemma}\label{lemmaVoronoiTiling}~{
    Разбиение $\mcT_{\mcA_dP}$ является разбиением на параллелоэдры Вороного.
  } \end{lemma}
  
  \smallskip
  \noindent {\it Доказательство.}
    Для стандартного параболоида $\Pi_I$ рассмотрим фокус $F$ и директрису $D$: точку 
    $\parenth{0, \ldots, 0, \frac14}^T \in \mathbb E^{d+1}$ и гиперплоскость $\braces{y^{d+1} = 
      -\frac14}$ соответственно. Известный факт, что $\Pi_I$ --- множество точек 
    равноудалённых от $F$ и $D$. Пусть точки $T_1$ и $T_2$ --- точки касания гиперплоскостей 
    $H_1$ и $H_2$ с $\Pi_I$, точка $M$ --- произвольная точка из $H_1\bigcap H_2$. $T''_1, T''_2, 
      M''$ --- ортогональные проекции соответствующих точек на директрису, а $T'_1, T'_2, M'$ --- 
    их ортогональные проекции на $\mathbb E^d = \braces{y^{d+1} = 0}$. По, так называемому, 
    оптическому свойству параболоида, $F$ и $T''_1$ симметричны относительно $H_1$, $F$ и $T''_2$
    симметричны относительно $H_2$. Так как $M \in H_1 \bigcap H_2$, то $MT''_1 = MF = MT''_2$.
    Отсюда, очевидно, $M''T''_1 = M''T''_2$. Из параллельности гиперплоскостей $\braces{y^{d+1} = 
      -\frac14}$ и $\braces{y^{d+1} = 0}$ следует, что и для проекций на $\mathbb E^d$ выполнено
    $M'T'_1 = M'T'_2$. Значит, в силу произвольности выбора $M$, проекция $H_1\bigcap H_2$ на 
    $\mathbb E^d$ --- это $(d-1)$-плоскость, равноудалённая от точек $T'_1$ и $T'_2$, то есть
    срединный перпендикуляр к отрезку $T'_1T'_2$.
    
    Пусть $T_1$ и $T_2$ --- точки касания произвольных двух ячеек $\widetilde{P_1}$ и 
    $\widetilde{P_2}$ из $\mcG_{\mcA_dP}$ с $\Pi_I$. Очевидно, $\widetilde{P_1}$ принадлежит
    одной из двух полуплоскостей, на которые $H_1\bigcap H_2$ разбивает $H_1$. Проекция этой 
    полуплоскости на $\mathbb E^d$ --- $d$-полуплоскость, ограниченная срединным перпендикуляром
    к отрезку $T'_1T'_2$: $\braces{x\in \mathbb E^d: \left\|x-T'_1\right\| \leqslant 
      \left\|x-T'_2\right\|}$. Пусть точка $T_1$ фиксирована, а $T_2$ пробегает все остальные точки
    касания женератрисы с $\Pi_I$. Пересечение построенных полуплоскостей гиперплоскости $H_1$
    равно самой грани $\widetilde{P_1} \subset \mcG_{\mcA_dP}$. Проекция данного 
    пересечения, с одной стороны --- это $\braces{x\in \mathbb E^d: \left\|x-T'_1\right\| \leqslant 
      \left\|x-\p\right\|, \forall \p\in \Lambda^d\parenth{\mcA_dP}}$ --- параллелоэдр Вороного для
    решётки $\Lambda^d\parenth{\mcA_dP}$. С другой стороны --- это проекция грани женератрисы, то 
    есть $d$-мерная ячейка разбиения $\mcT_{\mcA_dP}$.
  \hfill $\Box$ \par\bigskip
  
  \begin{theorem}\label{theoremMainEquivalence}~{ 
    Гипотеза Вороного выполнена для некоторого параллелоэдра $P$ тогда и только тогда, когда 
    существует каноническая нормировка разбиения $\mathcal T_P$.
  } \end{theorem}
  
  \smallskip
  \noindent {\it Доказательство.}
    Применяя последовательно лемму \ref{lemmaInvariantScaling}, теорему 
    \ref{theoremManifoldGeneratris}, теорему \ref{theoremTouch} и лемму \ref{lemmaVoronoiTiling} 
    получаем цепочку рассуждений: из существования некоторой канонической нормировки разбиения 
    $\mcT_P$ следует существование трансляционно инвариантной канонической нормировки. 
    Приведённая в разделе \ref{subsectionConstruction} конструкция по этой нормировке задаёт 
    корректно определённую женератрису, в которую вписан график соответствующей положительно 
    определённой квадратичной формы $Q$. Аффинное преобразование, приводящее график этой формы к 
    виду стандартного параболоида вращения, заданного единичной матрицей, переводит параллелоэдр 
    $P$ в параллелоэдр Вороного. В одну сторону теорема доказана.
    
    Для доказательства теоремы в обратную сторону отметим, что для произвольного разбиения на 
    параллелоэдры Вороного существует женератриса: достаточно поднять точки решётки центров
    $\Lambda^d(P_V)$ на стандартный параболоид $\Pi_I$ и провести через них касательные плоскости
    к параболоиду. Эти касательные плоскости ограничивают некоторый выпуклый $(d+1)$-полиэдр, 
    описанный около $\Pi_I$. Аналогично доказательству леммы \ref{lemmaVoronoiTiling} имеем, что
    проекция полученного полиэдра --- то же самое разбиение на параллелоэдры Вороного, $d$-мерная 
    граница полиэдра --- женератриса разбиения. Значит если нашлось аффинное преобразование, 
    переводящее параллелоэдр в параллелоэдр Вороного, то обратное аффинное преобразование даст 
    женератрису исходного разбиения. Отсюда, согласно теореме \ref{theoremManifoldGeneratris}, у 
    исходного разбиения есть каноническая нормировка.
  \hfill $\Box$ \par\bigskip
  
  \begin{corollary}[Теорема Житомирского]\label{corollaryZhitomirskii}~{
    Всякий $(d-2)$-примитивный параллелоэдр $P$ аффинно эквивалентен некоторому параллелоэдру 
    Вороного.
  } \end{corollary}
  
  \smallskip
  \noindent {\it Доказательство.}
    По лемме \ref{lemmaCSZhitomirskii} для разбиения $\mcT_P$ существует каноническая нормировка, 
    значит по теореме \ref{theoremMainEquivalence} для $P$ верна гипотеза Вороного.
  \hfill $\Box$ \par\bigskip

\subsection{Благодарности}

  Автор искренне благодарит Н.П. Долбилина, Р.М. Эрдала, А.И. Гарбера, А.Н. Магазинова за регулярные
  и плодотворные обсуждения и ценные замечания, В.П. Гришухина за редкие, но ценные разговоры и 
  предоставленные оттиски важных для написания данной работы статей Ч. Дэвиса и Ф. Ауренхаммера. 
  Данная работа была выполнена в замечательных для работы условиях в Queen's University в г. 
  Кингстон, Онтарио, Канада, коллективу которого автор также выражает свою искреннюю благодарность.

\end{document}